\documentclass{amsart}

\newtheorem{Proposition}{Proposition}
\newtheorem{Theorem}[Proposition]{Theorem}
\newtheorem{Lemma}[Proposition]{Lemma}

\usepackage{amssymb}

\newcommand{\notI}{\overline{I}}
\newcommand{\notJ}{\overline{J}}
\newcommand{\ssm}{\smallsetminus}

\def\max{{\textnormal{max}}}

\begin{document}

\bibliographystyle{../dart}

\title{A note on Bruhat order and double coset representatives}

\author{Christophe  Hohlweg}
\address[Christophe Hohlweg]{The Fields Institute\\
222 College Street\\
Toronto, Ontario, M5T 3J1\\ CANADA}
\email{chohlweg@fields.utoronto.ca}
\urladdr{http://www.fields.utoronto.ca/\~{}chohlweg}

\author{Mark Skandera}
\address[Mark Skandera]{Mathematics Department\\
Haverford College\\
370 Lancaster Ave\\
Haverford, PA 19041\\ USA}
\email{mskander@haverford.edu}
\urladdr{http://www.haverford.edu/math/skandera.html}

\maketitle

\section{Introduction}

Let $(W,S)$ be a finite Coxeter system with {\em length function}
$\ell$ and identity $e$. Endow
 $W$ with the {\em Bruhat order} $\leq$, that is, $w\leq g$ in $W$ if and only
  if an expression for $w$ can be obtained by deleting simple reflections in a reduced expression for $g$.
(If $w \leq g$ then we necessarily have $\ell(w)\leq \ell(g)$.)
    We refer the reader to
\cite{bjorner-brenti,humphreys} as general references for Coxeter
groups and the Bruhat order.

Denote by $W_I$ the {\em standard parabolic subgroup} of $W$
generated by $I\subset S$. For $I,J\subset S$, each double coset in
$W_I \backslash W / W_J$ has a unique minimal element.  Let 
$X_{IJ}=\{w\in W\,|\, w<rw,\ w<ws,\ \forall r\in I,\,\forall s\in J\}$
be the set of all minimal representatives of double
cosets in $W_I\backslash W/W_J$.

Curtis~\cite[Theorem~1.2]{curtis} shows that for any $I,J\subset
S$ and  $b\in X_{IJ}$, there is a unique maximal element  $b^\max$
 in $W_IbW_J$. This fact plays an important role in his study of Lusztig's isomorphism theorem.

The aim of this note is to prove the following result, which seems
to have escaped observation:

\begin{Theorem}\label{prop:Main} Let $I,J\subset S$ and $u,v\in X_{IJ}$.
Then $u\leq v$ if and only if $u^\max \leq v^\max$.
\end{Theorem}

Double parabolic cosets arise in a variety of settings.  In particular
 Theorem~\ref{prop:Main} is used in \cite{SkanNNDCB} in the study of the
dual canonical basis of $\mathcal{O}(SL_n\mathbb{C})$.\\

After proving our result in \S\ref{s2}, we give a combinatorial
criterion in \S\ref{s3} for the comparison
 of $u$ and $v$ (or $u^{\max}$ and $v^\max$.)

\section{Proof}\label{s2}

For $I\subset S$,  it is well-known
that the set
$$
  W^I =\{u\in W \,|\, u<us,\,\forall s\in I\}
$$
is a set of minimal length coset representatives of $W/W_I$.
 Each element  $w \in W$ has therefore a
unique decomposition $w = w^I w_I$ where $w^I \in W^I$ and $w_I
\in W_I$.
 Moreover
$\ell(w) = \ell(w^I) + \ell(w_I)$.  The pair $(w^I, w_I)$ is
generally refered to as the {\em parabolic components of $w$
along $I$} (see \cite[Proposition 2.4.4]{bjorner-brenti},
 or \cite[5.12]{humphreys}).  It is clear that $W^K_I=W^K\cap W_I$
 is a set of minimal length coset representatives of $W_I/W_K$. Moreover
 $X_{IJ}=(W^I)^{-1}\cap W^J$, where $(W^I)^{-1}=\{w^{-1}\,|\, w\in W^I\}$.

Let $w_{0,I}$ denote the unique maximal element in 
in $W_I$, and let $w_0 = w_{0,S}$ denote the longest element of $W$. 
Then the parabolic components of $w_0$ are $(w_0^I,w_{0,I})$, where 
$w_0^I$ is the unique maximal element
 in $W^I$ (see \cite[\S2.5]{bjorner-brenti}).
 It follows that for $K\subset I\subset S$, 
the unique maximal element in $W^K_I$
is $w_{0,I}^K=w_{0,I}w_{0,K}$. 

We recall the following well-known facts:

\begin{enumerate}

\item[(i)] For any $I,J\subset S$, define 
$I\cap bJb^{-1} = I \cap \{bsb^{-1}\,|\, s\in J\}$.  Then we have
\begin{equation*}\label{eq:DoubleCosetProd}
W^J = \coprod_{b\in X_{IJ}} W^{I\cap bJb^{-1}}_I b.
\end{equation*}
Therefore, each element $w\in W$ has a unique decomposition $(a,b,w_J)$ where
 $b\in X_{IJ}$, $a \in W^{I\cap bJb^{-1}}_I$
and $ab=w^J$.
Moreover $\ell(w)=\ell(a)+\ell(b)+\ell(w_J)$.
(See for instance \cite[\S2]{bbht}.)

\item[(ii)] Let $w,g,x\in W$ satisfy $\ell(wx)=\ell(w)+\ell(x)$
 and $\ell(gx)=\ell(g)+\ell(x)$. Then  $w\leq g$  if and only if  $wx\leq gx$.

\item[(iii)] From (ii) we have:
if $w\leq g$, then $w^I\leq g^I$ for any $w,g\in W$
and $I\subset S$.


\item[(iv)]  {\em Deodhar's Lemma  \cite{deodhar}}: Let $K\subset S$, $x\in W^K$ and $s\in S$. If $sx<x$
 then $sx\in W^K$. If $x<sx$ then either $sx\in W^K$ or  $sx=xr$ with $r\in K$.

\item[(v)] {\em Lifting property}:  Let $w,g\in W$ and $s\in S$  satisfy $w<sw$ and $sg<g$.
 Then $ w\leq g \iff w\leq sg \iff sw \leq g$ (see \cite{humphreys,bjorner-brenti}).
\end{enumerate}

Curtis~\cite[Theorem~1.2]{curtis} shows that for any $I,J\subset
S$ and  $b\in X_{IJ}$, $b^\max = w_{0,I}^{I\cap bJb^{-1}} b
w_{0,J}$ is the unique maximal
 element in $W_IbW_J$. Here we  give a short proof of this fact.   Let $w\in W_I bW_J$, then by (i) we have
$w=abw_J$ with
 $a\in W^{I\cap bJb^{-1}}_I$. Hence  $a\leq w^{I\cap bJB^{-1}}_{0,I}$
 and $w_J\leq w_{0,J}$, and by (i) and (ii) we have $w\leq b^\max$.

\begin{Lemma}\label{lem:1} Let $I,J\subset S$ and suppose that $u,v\in X_{IJ}$ 
satisfy $u\leq v$.  Then 
for any $a\in W_{I}^{I\cap uJu^{-1}}$ 
we have $au\leq w_{0,I}^{I\cap vJv^{-1}} v$.
\end{Lemma}

\begin{proof} Writing $\alpha=w_{0,I}^{I\cap vJv^{-1}}$, we will use induction
on $\ell(a)$ to show that $au \leq \alpha v$. 
If $\ell(a)=0$, then $a=e$ and $au=u\leq \alpha v$ since $u\leq v$ and $\ell(\alpha v)=
  \ell(\alpha)+\ell(v)$.

Assume therefore that $\ell(a)>0$. 
Then some $s\in I$ satisfies 
$sa<a$ and we have $sa\in W_I^{I\cap uJu^{-1}}$ by Deodhar's Lemma.
Since $\ell(sa)<\ell(a)$, we have
$sau\leq \alpha v$ by our induction hypothesis.
We also have $sau<au$ by (ii) since $a,sa\in W_I$ and $u^{-1}\in W^I$.  
In order to compare $\alpha v$ and $au$ we consider two cases.

If $s\alpha v < \alpha v$ then we obtain $au\leq \alpha v$ from (v) using 
$w=sau$ and $g=\alpha v$.
If $\alpha v<s \alpha v$, then $au=s(sau)\leq s\alpha v$ by definition. 
Observe that $\alpha<s\alpha$ by (ii). 
As $\alpha$ is the maximal element in $W_I^{I\cap vJv^{-1}}$, we have
$s\alpha \notin W_I^{I\cap vJv^{-1}}$ by (iv). So some $r\in I\cap vJv^{-1}$
satisfies $s\alpha = \alpha r$. 
Set $t=v^{-1}r v\in J$ so that $s\alpha v=\alpha v t$.
As $\alpha v\in W^J$ (by (i)) we deduce that $(\alpha v,t)$ 
are the parabolic components of
$s\alpha v$ along $J$.  As $au\in W^J$ and $au\leq s\alpha v$, 
we obtain by (iii) that $au=(au)^J \leq (s\alpha v )^J=\alpha v$.
\end{proof}

\begin{proof}[Proof of Theorem~\ref{prop:Main}]  
Assume that $u^\max\leq v^\max$. First observe that from (1)
and (i) we have 
$(((b^\max)^J)^{-1})^I=b^{-1}$ 
for any $b\in X_{IJ}$.
Now use (iii) and the automorphism $w \mapsto w^{-1}$ of the Bruhat 
order 
to show  that $u^\max\leq v^\max$ implies $u\leq v$.

Assume now that $u\leq v$. Using (1), (i) and (ii) we just have to
show that $w_{0,I}^{I\cap uJu^{-1}}u=(u^\max)^J\leq
(v^\max)^J=w_{0,I}^{I\cap vJv^{-1}}v$. But this is the 
special case $a=w_{0,I}^{I\cap uJu^{-1}}$ of Lemma~\ref{lem:1}.
\end{proof}

\section{The special case of the symmetric group}\label{s3}

Given a permutation $w \in S_n$, we define the matrix $M(w) = (m_{i,j}(w))$ by
setting $m_{i,j}(w) = \delta_{j,w_i}$, where $w_1 \cdots w_n$ is the one-line
notation of $w$.  We define the related matrix
$D(w) = (d_{i,j}(w))$ by
$d_{i,j}(w) = \sum_{k = 1}^i \sum_{\ell = 1}^j m_{i,j}(w)$.
It is well known that $u \leq v$ if and only if we have the componentwise
inequality of matrices $D(u) \geq D(v)$, and we shall state a similar fact
for double parabolic analogs of $M$ and $D$.

A subset $H$ of $S = \{ s_1,\dotsc, s_{n-1} \}$
induces an equivalence relation $\sim_H$
on $[n] = \{1,\dotsc,n\}$ which is the transitive closure of the relation
$i~R~(i+1)$ for all $s_i \in H$. Let $B_1, \dotsc, B_p$ and $C_1 ,\dotsc, C_q$
be the equivalence classes of $\sim_I$ and $\sim_J$, respectively,
and define the matrices
$M^{I,J}(w) = (m_{i,j}^{I,J}(w))$ and $D^{I,J}(w) = (d_{i,j}^{I,J}(w))$ by
\begin{equation*}
m_{i,j}^{I,J}(w) = \# \{ k \in B_i \,|\, w_k \in C_j \},
\qquad
d_{i,j}^{I,J}(w) = \sum_{k=1}^i \sum_{\ell=1}^j m_{i,j}^{I,J}(w).
\end{equation*}
It is well known (see, e.g., \cite{JK}) that $u$ and $v$ belong to the same
double coset in $W_I \backslash W / W_J$
if and only if $M^{I,J}(u) = M^{I,J}(v)$.  Furthermore we have the following.

\begin{Proposition}\label{prop:matrix}
Given $u,v$ in $X_{I,J}$, then $u \leq v$ (or $u^\max \leq v^\max$)
if and only if we have the componentwise inequality of matrices
$D^{I,J}(u) \geq D^{I,J}(v)$.
\end{Proposition}
\begin{proof}
Define $\notI = [n] \ssm \{ i \,|\, s_i \in I \}$,
$\notJ = [n] \ssm \{ j \,|\, s_j \in J \}$.
Then for each $w \in S_n$, the matrix $D^{I,J}(w)$ is equal to the
$(\notI, \notJ)$ submatrix of $D(w)$.  The ``only if'' direction follows
immediately.

Suppose the $u \not \leq v$ and let $(i,j)$ be a componentwise minimal
pair satisfying $d_{i,j}(u) < d_{i,j}(v)$.
If $i > 1$, then the fact that the matrices $D(u)$ and $D(v)$
weakly increase down columns and across rows,
with adjacent entries differing by no more than $1$, implies that
$d_{i-1,j}(u) \leq d_{i,j}(u) < d_{i,j}(v) \leq d_{i-1,j}(v) + 1$.
By the minimality of $i$ and $j$, this last expression is less than
or equal to $d_{i-1,j}(u) + 1$, and for some nonnegative integer $c$ we have
$d_{i-1,j}(u) = d_{i,j}(u) = d_{i-1,j}(v) = c$
and
$d_{i,j}(v) = c+1$.
Similarly, if $j > 1$ then we have
$d_{i,j-1}(u) = d_{i,j}(u) = d_{i,j-1}(v) = c$.
It follows that for any values of $(i,j)$ we must have
$u_i > j$ and $u_j^{-1} > i$.

Now let $(k,\ell)$ be the componentwise minimal pair in $\notI \times \notJ$
satisfying $i \leq k$, $j \leq \ell$.
Since $u \in X_{I,J}$, we must also have
\begin{equation*}
\ell < u_i < \cdots < u_k,\qquad
k < u_j^{-1} < \cdots < u_\ell^{-1}.
\end{equation*}
Thus $d_{k,\ell}(u) = c$.
Since $d_{k,\ell}(v) \geq c+1$, we conclude that
$D^{I,J}(u) \not \geq D^{I,J}(v)$.
The equivalence of $D^{I,J}(u) \geq D^{I,J}(v)$ and $u^\max \leq v^\max$ 
follows from a similar argument.
\end{proof}

We illustrate Proposition~\ref{prop:matrix} by considering
$W = S_7$,
subsets
$I=\{s_1,s_2,s_4,s_6 \}$, $J=\{s_1,s_3,s_4,s_5\}$
of generators, 
and
corresponding 
equivalence classes $123|45|67$, $12|3456|7$.
To compare minimal 
representatives $u = 1342567$,
$v = 3471526$ of two double cosets in $W_I \backslash W / W_I$,
we use the matrices
\begin{equation*}
M^{I,J}(u) = \begin{bmatrix}
1 & 2 & 0 \\
1 & 1 & 0 \\
0 & 1 & 1
\end{bmatrix},
\qquad
M^{I,J}(v) = \begin{bmatrix}
0 & 2 & 1 \\
1 & 1 & 0 \\
1 & 1 & 0
\end{bmatrix},
\end{equation*}
to compute 
\begin{equation*}
D^{I,J}(u) = \begin{bmatrix}
1 & 3 & 3 \\
2 & 5 & 5 \\
2 & 6 & 7
\end{bmatrix},
\qquad
D^{I,J}(v) = \begin{bmatrix}
0 & 2 & 3 \\
1 & 4 & 5 \\
2 & 6 & 7
\end{bmatrix},
\end{equation*}
and conlcude that $u \leq v$ and $u^\max \leq v^\max$.



\end{document}